\documentclass[reqno]{amsart}

\usepackage{enumerate}
\usepackage{ifpdf}
\usepackage{amsmath}
\usepackage{amsfonts}
\usepackage{amssymb}
\usepackage{amsthm}
\usepackage{hyperref}
\usepackage{setspace}
\usepackage{amsrefs}

\newcommand{\ignore}[1]{}

\renewcommand{\Re}{\operatorname{Re}}
\renewcommand{\Im}{\operatorname{Im}}

\newcommand{\abs}[1]{\left\lvert {#1} \right\rvert}

\newcommand{\C}{{\mathbb{C}}}
\newcommand{\R}{{\mathbb{R}}}

\newcommand{\sH}{{\mathcal{H}}}

\newcommand{\sX}{{\mathcal{X}}}

\newtheorem{thm}{Theorem}[section]
\newtheorem*{thmnonum}{Theorem}

\newtheorem{cor}[thm]{Corollary}
\newtheorem{lemma}[thm]{Lemma}

\theoremstyle{definition}

\theoremstyle{remark}

\theoremstyle{remark}
\newtheorem{example}[thm]{Example}

\author{Ji\v{r}\'i Lebl}
\address{ Department of Mathematics, University of California
at San Diego, La Jolla, CA 92093-0112, USA}
\email{jlebl@math.ucsd.edu}
\date{May 7, 2007}

\ifpdf
\fi

\title[Extension of Levi-flat hypersurfaces past CR boundaries]
{Extension of Levi-flat hypersurfaces past CR boundaries}

\begin{document}


\begin{abstract}
Local conditions on boundaries of $C^\infty$ Levi-flat hypersurfaces, in case
the boundary is a generic submanifold, are studied.  For nontrivial real
analytic boundaries we get an extension and uniqueness result, which forces the
hypersurface to be real analytic.  This allows us to classify all real
analytic generic boundaries of Levi-flat hypersurfaces in terms of their
normal coordinates.  For the remaining case of generic real analytic boundary
we get a weaker extension theorem.  We find examples to show that these
two extension results are optimal.
Further, a class of nowhere minimal real analytic submanifolds is found, which
is never the boundary of even a $C^2$ Levi-flat hypersurface.
\end{abstract}

\maketitle


\section{Introduction} \label{section:intro}

The question we wish to ask is when is a generic codimension 2 submanifold
$M \subset \C^N$ locally the boundary of a Levi-flat hypersurface $H$.
In particular, we will
ask the following questions.  When does $H$ extend as a Levi-flat
hypersurface past $M$?
When is $H$ unique?  How does the regularity of $H$
depend on the regularity of $M$?  We will answer these questions
fully when $M$ is real analytic and $H$ is smooth.

The results here are motivated by Dolbeault, Tomassini and Zaitsev
\cite{DTZ:boundfull}, who consider the
global situation under additional assumptions on $M$.
These results are also related to results of Straube and Sucheston
\cite{StraubeSucheston:fol}.
This paper can also
be seen as a natural extension of the results in \cite{Lebl:lfnm},
as we will mostly concern
ourselves with the situation when $M$ is real analytic.
In the non-CR case, which is not considered here, similar questions are
considered in $\C^2$ for example by Bishop \cite{Bishop:diffman},
Moser and Webster \cite{MW:normal}, or
Bedford and Gaveau \cite{BG:envhol}.  For further discussion of the
non-CR case and more references see \cite{DTZ:boundfull}.

In the following, by \emph{submanifold} we always mean embedded submanifold,
by \emph{hypersurface} a submanifold of codimension 1, and by
\emph{real analytic subvariety} of an open set $U$,
a set closed in $U$ and locally defined by
the vanishing of a family of real analytic functions.

Let $M \subset \C^N$ be a real codimension 2 connected submanifold.
Let $J$ be the complex structure on $\C^N$, and let
$T_p^cM = J(T_pM) \cap T_pM$.
A real submanifold is called \emph{CR} if the dimension of
$T_p^cM$ is constant as $p$ varies in $M$.
The smallest germ (in terms of dimension) of a
CR submanifold $N$ of $M$ through $p$ such that the $T_q^c N = T_q^c M$
is called the \emph{local CR orbit} at $p$, and is guaranteed to exist
by the Nagano theorem \cite{Nagano}
in case $M$ is real analytic, or the Sussmann theorem \cite{Sussmann} if
$M$ is only smooth.  $M$ is said to be \emph{minimal at $p$} (in the sense of
Tumanov \cite{Tumanov}) if the local
CR orbit through $p$ is of codimension 0 in $M$.  If $M$ is not minimal at
any point then $M$ is said to be \emph{nowhere minimal}.
See \cites{BER:book, Boggess:CR, DAngelo:CR} for more details.

We will say that $M$ is \emph{generic} if and only if
$T_pM + J(T_pM) = T_p\C^N$ for all $p \in M$, where $J$ is the complex
structure on $\C^N$.
If $M$ is a real analytic CR submanifold of
codimension 2, this just means that $M$ is not a complex analytic submanifold
near any point.
We will always assume that $0 \in M$.

A set $H \subset \C^N$ is
a \emph{$C^k$ hypersurface with boundary}, if
there is a subset $\partial H \subset H$,
such that
$\partial H \subset H$,
$H \setminus \partial H$ is a $C^k$ hypersurface
(submanifold of codimension 1), and for each point $p \in
\partial H$,
there exists a neighbourhood $p\in U \subset \C^N$,
a $C^k$ diffeomorphism $\varphi \colon U \to \R^{2N}$,
such that
$\varphi (H \cap U) = \{ x \in \R^{2N} \mid x_{2N-1} \geq 0, ~ x_{2N} = 0\}$,
and such that
$\varphi (\partial H \cap U) = \{ x \in \R^{2N} \mid x_{2N-1} = 0,~ x_{2N} =
0\}$.  Hence, $\partial H$ is a $C^k$ submanifold of codimension 2
in $\C^N$.
We will call $H^o := H \setminus \partial H$ the interior of $H$.
As we are concerned
with only local questions, we can assume that there exists just one such
$U$, and such that $\partial H, H \subset U$.  We can further assume that
$\partial H$ and $H$ are closed subsets of $U$.
We can extend $H$ to $\tilde{H}$, a full $C^k$ submanifold
near $0$, by just pulling back a neighbourhood of $0 \in \R^{2N}$ by
$\varphi$.

A $C^k$ ($k \geq 2$) hypersurface $H$ is said to be \emph{Levi-flat}
if the bundle $T^cH$ is
involutive\footnote{This is equivalent to saying that the Levi form vanishes
identically at every point, which is the usual definition.}.
An equivalent definition is to say that near every
point of $H$, there exists a one parameter
local foliation of $H$ by complex hypersurfaces, which is called the
\emph{Levi foliation}.  To see why these are equivalent, note that
if $T^cH$ is involutive
the Frobenius theorem gives us a $C^{k-1}$
foliation with the leaves being complex hypersurfaces (they are locally the
graphs of holomorphic functions).
If $H$ is a hypersurface with boundary as defined above, then 
we will say it is Levi-flat when $H^o$ is Levi-flat.
If $H$ is a real analytic
subvariety of codimension 1, then we say it is Levi-flat, if it is Levi-flat as a submanifold
at all the nonsingular points.
We can now state our main result.

\begin{thm} \label{extthm}
Let $M \subset \C^N$
be a connected real analytic generic submanifold of codimension 2
through the
origin, such that not all local CR orbits of $M$ are of codimension 2 in $M$.
Suppose that there exists a
connected Levi-flat $C^\infty$ hypersurface $H$ with boundary, 
where
$M \subset \partial H$.
Then there exists
a neighbourhood $U$ of the origin and a nonsingular real analytic
Levi-flat hypersurface $\sH$ such that $H \cap U \subset \sH$.

Further, the germ $(\sH,0)$ is unique in the sense that if
$(\sH',0)$ is a germ of an
irreducible
real analytic Levi-flat subvariety of codimension 1
such that $(M,0) \subset (\sH',0)$,
then $(\sH',0) = (\sH,0)$.
\end{thm}

First, note that the condition that $M$ is real analytic is necessary for
the extension to hold.
See Example \ref{example:smooth}
for a counterexample in case $M$ is $C^\infty$.

The condition on the local CR orbits is necessary for the conclusion that
the extension
$\sH$ is unique and real analytic.  If $M$ is the boundary of a Levi-flat
hypersurface, then all local CR orbits must be of positive codimension in $M$,
see Lemma \ref{lemma:nomin}.  If all the local CR orbits are of codimension
1 in $M$, then the theorem follows easily by known results,
see Lemma \ref{lemma:uniqueH}.
Finally, if all local CR orbits would be of codimension
2 in $M$, then $M \subset \C^N$ would be locally biholomorphic to
$\C^{N-2} \times \R^2$,
and we will give (Example \ref{example:smoothH})
an example of a bona fide
$C^\infty$ (i.e. not contained in a real analytic subvariety)
Levi-flat hypersurface which contains such an $M$.
Hence the theorem is, in this respect, optimal.
In \S \ref{section:flatex}, we will prove the following weaker
extension theorem for such submanifolds, which is also optimal in view of the
above examples.
In the sequel, when we consider $\C^{N-2} \times \R^2$ as a subset
of $\C^N$, we mean the natural embedding.

\begin{thm} \label{flatextthm}
Suppose $H \subset \C^N$ is a $C^\infty$ Levi-flat
hypersurface with boundary, and $0 \in \partial H \subset \C^{N-2} \times \R^2$.
Then for some neighbourhood $U$ of the origin,
there exists a $C^\infty$ Levi-flat hypersurface $\sH$ (without boundary)
such that $H \cap U \subset \sH$.

Further, the germ $(\sH,0)$ is unique in the sense that if
$(\sH',0)$ is another a germ of a $C^\infty$
Levi-flat hypersurface such that $(H,0) \subset (\sH',0)$,
then $(\sH',0) = (\sH,0)$.
\end{thm}

Note that the uniqueness in Theorem \ref{flatextthm}
is much weaker as $\sH$ depends on $H$,
whereas in Theorem \ref{extthm} $\sH$ depends only on $M$.

Theorem \ref{extthm} says that
in particular,
there exists a holomorphic function defined near the origin
with nonzero gradient that is real valued on $M$.  In other
words, $M$ is locally the
boundary of a Levi-flat $C^\infty$ hypersurface if and only if $M$
has local defining functions in $(z,w) \in \C^{N-2} \times \C^2$ of the form:
\begin{equation} \label{eq:nicecoord}
\begin{split}
& \Im w_1 = \varphi(z,\bar{z},\Re w) ,\\
& \Im w_2 = 0 ,
\end{split}
\end{equation}
for some $\varphi$ such that
$\varphi(0,\bar{z},s) \equiv \varphi(z,0,s) \equiv 0$
(i.e.\@ these are normal coordinates, see \cite{BER:book} for example).
The classification of Levi-flat boundaries that are generic and
real analytic is therefore simple.

\begin{cor}
Let $M \subset \C^N$
be a connected real analytic generic submanifold of codimension 2
through the origin.  The following are equivalent:
\begin{enumerate}[(i)]
\item
There exists a Levi-flat $C^\infty$
hypersurface $H$ with boundary, such that $0 \in \partial H \subset M$.
\item
There exists a real analytic Levi-flat hypersurface (submanifold) $H$
defined in a neighbourhood $U$ of the origin such that
$M \cap U \subset H$.
\item
There exist local holomorphic coordinates (near
the origin) such that $M$ is defined by an equation of the form
\eqref{eq:nicecoord}.
\item
There exists a real analytic foliation of codimension 1 in $M$,
defined in a neighbourhood of the origin, such that the leaves are
unions of (representatives of) local CR orbits of $M$.
\end{enumerate}
\end{cor}

When the Levi-flat hypersurface is only $C^2$ rather than smooth, then we
will be able to prove that the individual leaves of the Levi foliation
extend across $M$.
See Lemma \ref{lemma:leavesextend}.
As an application of this lemma we prove the
following theorem.  First,
we must define the property of being {\em almost minimal} (see
\cite{Lebl:lfnm}).
Let $M$ be a
real analytic,
generic submanifold through the origin.  Suppose that for every $U$
a neighbourhood of the origin there exists a $p \in U \cap M$ such that the
local CR orbit (take a representative of this germ) of $M$ through $p$ is not
contained in any complex analytic
subvariety of $U$, then $M$ is almost minimal at the origin.
An example of this kind of manifold can be found in
\S \ref{section:ex}.

\begin{thm} \label{noalmostmin}
Let $M \subset \C^N$ be a connected real analytic generic submanifold of
codimension 2 through the origin, which is almost minimal at the origin.
Let $H$ be a connected $C^2$ hypersurface with boundary and $M \subset
\partial H$.  Then $H$ is not Levi-flat.
\end{thm}

If $H$ would be $C^\infty$ then the above result follows at once from
Theorem \ref{extthm}.  Further, not being almost minimal is a
necessary, but not sufficient, condition to being a boundary of a $C^2$
Levi-flat hypersurface.

The organization of this paper is as follows.  In \S \ref{section:flatbnd}
we discuss boundaries of Levi-flat hypersurfaces in general and prove 
Theorem \ref{extthm}.  In \S \ref{section:flatex} we prove
Theorem \ref{flatextthm}.  In \S \ref{section:counterex}
we give examples that show that Theorems \ref{extthm} and
\ref{flatextthm} are optimal.
In \S \ref{section:almostmin} we prove Theorem \ref{noalmostmin}.
In \S \ref{section:ex} we give an example almost minimal submanifold
which does not ``bound'' (in a very weak sense) even a singular
Levi-flat real analytic
subvariety.  Finally, in \S \ref{section:subanal} we discuss the existence
of subanalytic Levi-flat hypersurfaces.

The author would like to acknowledge Dmitri Zaitsev for many useful comments
and suggestions to simplify certain proofs and improve exposition.
The author would also like to
acknowledge Peter Ebenfelt for many useful discussions and many suggestions
for improvements to these results.
Finally the author would like to thank the referee for pointing out a problem
with the statement and proof of Theorem \ref{C2extthm}, as well as other
comments and suggestions.

\section{Locally flat boundaries} \label{section:flatbnd}

We prove some basic results about locally flat boundaries.
For the rest of this section, we assume that $H$ is a hypersurface
with boundary, that $M = \partial H$,
and
that $M$ is a generic submanifold through the origin.

\begin{lemma} \label{lemma:nomin}
Let $M$ be $C^\infty$ and $H$ be $C^2$, and suppose 
that $H$ is Levi-flat, then $M$ is nowhere minimal.
\end{lemma}

\begin{proof}
We can just extend $H$ to $\tilde{H}$ as in the introduction and assume $\rho$
is a defining function for $\tilde{H}$.  Then $\theta = i(\partial \rho -
\bar{\partial} \rho)$ is a real $C^1$ one-form
that vanishes on $T^cH$.  On $H$, as $H$ is Levi-flat,
$d\theta \wedge \theta = 0$ and by continuity this happens on $M$ as well
if we restrict $\theta$ to $M$.  $\theta$ can't vanish on $M$
as that would make $M$ have a complex tangency (it would be tangent to
$T^c\tilde{H}$).  Hence there exists (locally at near every
point) a foliation of $M$ by CR
submanifolds of smaller dimension with the same CR dimension as $M$, and so
$M$ cannot be minimal at any point.
\end{proof}

\begin{lemma} \label{lemma:foliation}
Let $M$ and $H$ be $C^k$ ($2 \leq k \leq \infty$), and suppose 
that $H$ is Levi-flat, then the Levi foliation of $H^o$ extends to 
a foliation of $H$.  That is, in a perhaps a smaller neighbourhood of the
origin, there exists a $C^{k-1}$, real valued, function $f$ on $H$
(including $M$)
with nonvanishing differential ($f|_M$ also has nonvanishing differential),
such that $f$ is constant along leaves of
the Levi foliation of $H^o$.  If $M$ and $H$ are $C^\infty$, then $f$ is
$C^\infty$.
\end{lemma}

\begin{proof}
If $k=\infty$,
then by $C^{k-1}$ we will mean $C^\infty$ below.
For convenience we change notation slightly.
We straighten out the boundary, and assume $H$ is the upper half
plane $\{ x \in \R^n \mid x_1 \geq 0 \}$ and $M$ is defined by $x_1 = 0$
(where $n=2N-1$).
The $C^{k-1}$ 1-form given in the proof
of Lemma \ref{lemma:nomin} that vanishes on the vectors in $T_p^c H$
induces a
$C^{k-1}$ 1-form
$\theta$
on the upper half plane in $\R^n$.  $\theta$ does not vanish on the tangent
vectors to $x_1 = 0$
(else $M$ would have a complex tangency).  We can easily extend $\theta$
to all of $\R^n$ (or at least a neighbourhood of the origin) 
as a $C^{k-1}$ 1-form.
We now follow the proof of the Frobenius theorem in \cite{Flanders:dforms},
to show that there exists a real valued function with nonvanishing
differential at 0 that is constant on the Levi foliation of $H^o$.
That is, we just need to show that we can modify $\theta$ on the set
$x_1 < 0$, such that the modification is completely integrable.
We have that $d\theta \wedge \theta = 0$ for $x_1 \geq 0$.
It is not hard to see that there exists a $C^{k-2}$ 1-form $\alpha$ defined
near the origin such that $d \theta = \theta \wedge \alpha$ for $x_1 \geq 0$.

As $\theta$ does not vanish near the origin (and does not vanish identically
on $T_0M$),
we may assume that $\theta = dx_n + \sum_{j=1}^{n-1} A_j dx_j$.  Fix a point
$a$ in $x'$ space, where $x' = (x_1,\ldots,x_{n-1})$.  We consider the 
equation $\theta = 0$ on the hyperplane where $x_j = a_j t$ for $t \in \R$.  We
solve this ODE for $x_n$, with the initial condition $x_n(0) = c$, for some
constant $c$.  That is, we find the unique solution of
\begin{equation}
\begin{split}
& \frac{\partial F}{\partial t}(t,a,c) = - \sum_j A_j(at,F(t,a,c)) a_j , \\
& F(0,a,c) = c .
\end{split}
\end{equation}
We note that we can change scale $F(t,a,c) = F(kt,a/k,c)$, and hence setting
$k=1/t$, we get $F(t,a,c) = F(1,ta,c)$.  We change variables to $(u,v) \in
\R^{n-1} \times \R$ by
\begin{equation}
\begin{split}
& x' = u , \\
& x_n = F(1,u,v) .
\end{split}
\end{equation}
It is not hard to check that this is a change of coordinates.
In these new coordinates we write
\begin{equation}
\theta = \sum_j P_j du_j + B dv .
\end{equation}
Now we define
\begin{equation}
\tilde{\theta} = B dv .
\end{equation}
If we show that the $P_j$
vanish for $x_1 \geq 0$ ($u_1 \geq 0$), then we are done.
We know that $\sum P_j(ta,v)a_j = 0$.  This implies that if we consider
the mapping $\varphi(t,a,v) := (ta,v)$, we get
\begin{equation}
\varphi^* \theta = \sum \tilde{P}_j (t,a,v) da_j + \tilde{B}(t,a,v) dv .
\end{equation}
In particular, $\varphi^* \theta$ does not depend on $dt$.  Further
$\tilde{P}_j(t,a,v) = tP_j(ta,v)$ so $\tilde{P}_j(0,a,v) = 0$.
Now suppose that $a_1 \geq 0$ and $t \geq 0$, then
we have that $d(\varphi^*\theta) = (\varphi^* \alpha) \wedge
(\varphi^* \theta)$.  We set $D$ so that
$\varphi^* \alpha = D(t,a,v)dt + \ldots$.  From this equation we obtain
\begin{equation}
\frac{\partial \tilde{P}_j}{\partial t} = D \tilde{P}_j .
\end{equation}
By the uniqueness theorem for ODEs and the fact that
$\tilde{P}_j(0,a,v) = 0$ this implies that $\tilde{P}_j$ is identically zero,
and hence $P_j$ is identically zero.  This was true for $a_1 \geq 0$ ($t \geq
0$) and hence on the upper half plane and hence on $H$.  We therefore
have $\theta = \tilde{\theta}$ on the upper half plane and $\tilde{\theta}$
is closed and thus exact.  We get our $f$ of class $C^{k-1}$ (or
$C^\infty$ if $k=\infty$)
by taking $v$ as a function of $x$.
\end{proof}

\begin{lemma} \label{lemma:uniqueH}
Let $M$ be real analytic and $H$ be $C^2$, and
suppose that the local CR orbits of $M$ are all of codimension 1 in $M$.
Then there exists a neighbourhood $U$ of the origin such that
$(U \cap H) \subset \sH$, where $\sH$ is the unique Levi-flat real analytic
hypersurface in $U$ that contains $M$.
\end{lemma}

Note that $\sH$ is the union of the intrinsic complexifications of the local
CR orbits of $M$.  Where the intrinsic complexification is the smallest
complex submanifold containing the local CR orbit.

\begin{proof}
Since $M$ is real analytic and the local CR orbits are all of codimension 1
in $M$,
we can therefore apply the analytic Frobenius theorem to get a real analytic
real valued function on some small neighbourhood $U$ of the origin in $M$
with nonvanishing differential that is constant along the local CR orbits of
$M$.  Such a function is CR and hence extends to be holomorphic and the
vanishing of its imaginary part defines a Levi-flat hypersurface $\sH$.

Assume that $H \subset U$.
We must show that $H \subset \sH$.  By Lemma \ref{lemma:foliation},
we have that the Levi foliation of $H^o$ extends to $M$ (by perhaps making
$U$ smaller still).  That is, we have
complex submanifolds of $\C^n$ with boundary on $M$.  It is not hard to see
by the arguments used above that a leaf $L \subset H$ extended to the
boundary intersects $M$ precisely on a local CR orbit (by dimension).
The function that defines the corresponding leaf of the Levi foliation of
$\sH$
is of course holomorphic on $L$ and zero on the boundary of $L$, hence
$L \subset \sH$, and so $H \subset \sH$.
\end{proof}

\begin{proof}[Proof of Theorem \ref{extthm}]
Let $L_k$, $k = 1,\ldots,2N-4$,
be a basis of real analytic vectorfields spanning $T_p^cM$ defined
near the origin.  As  not all local CR orbits are of codimension 2 in $M$, then
there must exist an iterated commutator $K$ of the $L_k$, which is not
identically zero.  As $M$ is nowhere minimal (by Lemma
\ref{lemma:nomin}), then by dimension, $K$ together
with $L_k$ span the tangent space of the CR orbit whenever $K$ is nonzero.

By Lemma \ref{lemma:foliation} we have a $C^\infty$ codimension 1
foliation on $M$.  Hence, by forgetting for a moment the CR structure
of $M$, we can reduce to a situation
where we have a $C^\infty$ codimension 1
foliation on a small
neighbourhood $U \subset \R^{2N-2}$, given by a $C^\infty$
submersion $\varphi \colon U \to \R$, and
real analytic vector fields $L_k$ and $K$,
which are tangent
to the leaves of the foliation, $L_k$ never vanish and $K$ does not vanish
identically.
To see that the foliation must be real analytic, we
only need
to look at $TU$, the tangent bundle of $U$, and look at the normal bundle of
the foliation:
\begin{equation}
\{ (x,v) \in U \times \R^{2N-2} = TU \mid
\nabla \varphi (x) = t v, t \in \R \} ,
\end{equation}
which is a $C^\infty$ submanifold of dimension $2N-1$.
We define a larger real analytic subvariety of the same dimension:
\begin{equation}
\{ (x,v) \in U \times \R^{2N-2} = TU \mid K(x) \cdot v = 0 , L_k(x) \cdot v =
0, k = 1,\ldots,2N-4 \} ,
\end{equation}
where we view $L_k$ and $K$ as an $\R^{2N-2}$ valued function, and
the dot is the usual dot product.
Hence by a theorem 
of Malgrange (see \cite{Malgrange} Chapter VI, Proposition 3.11),
we see that the normal bundle to the foliation must be a real analytic
submanifold.  Therefore there must exist (locally near the origin, by Frobenius)
a real valued, real analytic submersion $f \colon M \to \R$
defining the foliation.  This submersion is constant along the local CR
orbits of $M$ and hence must be a CR function.  All real analytic CR functions
extend
uniquely to holomorphic functions in $\C^N$.  Thus $f$ is really a
holomorphic function with a nonvanishing gradient on $M$, which is real valued
on $M$.  Hence the equation $\Im f = 0$ defines a real analytic
Levi-flat hypersurface $\sH$, which contains $M$.
$\sH$ must contain
$H$ since it must contain the leaves of the Levi foliation of
$H$ by Lemma \ref{lemma:uniqueH},
and the leaves of $H$ are given by the foliation given by Lemma
\ref{lemma:foliation}.  Actually, Lemma \ref{lemma:uniqueH} only tells us
about leaves that pass through points of $M$ where the codimension in $M$ of
the local CR orbit is 1.  However, the remaining points lie on a real analytic
subvariety of $M$, and hence leaves that only pass through these points are
isolated and thus must also lie in $\sH$, since it is locally closed.

The uniqueness of $\sH$ is one of the conclusions of Theorem 1.1
in \cite{Lebl:lfnm}.
\end{proof}

\section{Extension across flat boundaries} \label{section:flatex}

When the local CR orbits of $M$ are all of codimension 2 in $M$, the situation
is different.  In this section we will prove Theorem \ref{flatextthm}.
First we will prove this result in $\C^2$, and then reduce the general
case to this.
In \S \ref{section:counterex}, we will see that a $C^\infty$ extension
is the best we can do.
Suppose that $\tau$ is the complex conjugation function.

\begin{thm} \label{C2extthm}
Suppose that
$H \subset \C^2$
is a Levi-flat $C^\infty$
hypersurface with boundary, with $0 \in \partial H \subset \R^2$.
Then there exists a neighbourhood
$U \subset \C^2$ of $0$, with $U = \tau (U)$ such that
$(H \cap U) \cup \tau(H \cap U)$ is a
$C^\infty$ Levi-flat hypersurface (without boundary).
\end{thm}

The idea is to extend the leaves of the Levi foliation of $H$ across $\R^2$.
Because $H$ has a boundary on $\R^2$, the leaves must be subvarieties of
$U \setminus \R^2$, and further can be extended to be complex submanifolds
of $U$.

\begin{proof}
Let $(z,w) \in \C^2$ be the coordinates.  As in the introduction,
let $\tilde{H}$ be any $C^\infty$
hypersurface without boundary through $0$, such that $H \subset \tilde{H}$.
Since $\partial H \subset \R^2$,
then either the $z$ or $w$ imaginary axis is not tangent to $\tilde{H}$,
so let us assume it is the $w$ imaginary axis.
Then there exists a
polydisc $U$ with center at $0$
such that $\tilde{H} \cap U$ is graph of a real valued continuous function
over $\{ \Im w = 0 \} \cap U$ (i.e. a function in $z$ and $\Re w$). 
Further we can choose $U$ such that $\tilde{H} \cap U$ and $H \cap U$
are closed in $U$.  To simplify notation, assume
$\tilde{H} = \tilde{H} \cap U$ and
$H = H \cap U$.
Now $\R^2 \cap U$ is a subset of $\tilde{H}$, and thus $H$ is without loss of
generality a graph
over $\{ \Im w = 0 \text{ and } \Im z \geq 0 \} \cap U$.
This means that $\tau(H) \cap H \subset \R^2$.
Hence
$\tau(H) \cup H$ is a graph of a continuous real function over
$\{ \Im w = 0 \} \cap U$.  Thus it remains to be shown that near every point
$p \in \partial H$, the union $\tau(H) \cup H$
is a smooth submanifold.

By Lemma \ref{lemma:foliation} the foliation of $H$ extends up to $\R^2$.  In
particular the leaves (extended to the boundary) are closed subsets of $H$.
Let $L \subset U$ be a leaf of the foliation of $H$ extended to the
boundary of $H$ as a submanifold with boundary.
Pick $p = (z_0,w_0) \in \partial H \subset \R^2$.
Since $L$ is a submanifold with boundary, we look at any $\tilde{L}$ being
an arbitrary $C^\infty$ extension of $L$ and we can assume $\tilde{L} \subset
\tilde{H}$.  We know from before that $\tilde{L}$ is not tangent to the
$w$ imaginary axis at $p$ as $\tilde{H}$ is not.  Since $T^c_p\tilde{L} =
T_p\tilde{L}$, we know that $\tilde{L}$ is not tangent to the real $w$ axis
near $p$ either.  Hence there exists a small polydisc 
$V$ with center $(z_0,w_0)$,
such that $\tilde{L} \cap V$ is a graph
of a continuous complex valued function of $z$ over $\{ w = w_0 \} \cap V$.
Since
$L\cap V$ lies above $\{ w = w_0 \text{ and } \Im z \geq 0 \} \cap V$
because $L \subset H$, then
as before $\tau(L \cap V)$ is a graph over
$\{ w = w_0 \text{ and } \Im z \leq 0 \} \cap V$ and
$(L \cup \tau(L)) \cap V$ is a graph of a continuous complex valued
function of $z$ over $\{ w = w_0 \}$.
Further this function is real
valued when $\Im z = 0$, and holomorphic when $\Im z > 0$.  Hence by Schwarz
reflection principle it is holomorphic everywhere and
$(L \cup \tau(L)) \cap V$ is a complex submanifold.  Since this is true
near all $p \in \R^2$, then $L \cup \tau(L)$
is a complex submanifold.

We can therefore foliate the set $H \cup \tau(H)$ by the complex submanifolds
$L \cup \tau(L)$.
Since the leaves of the foliation are complex submanifolds (and hence
$C^\infty$) and are not tangent to $\R^2$, then
$H \cup \tau(H)$ must be a $C^\infty$ submanifold.
To see this, we recall that
$H \cup \tau(H)$ is a graph of a real function $f$ over 
$\{ \Im w = 0 \} \cap U$ and let $\pi$ be the projection onto this plane.
Further, note that $H$ and $\tau(H)$ are $C^\infty$ up to the boundary
and hence all partial derivatives of $f|_{\pi(H)}$ extend
to $\R^2$ and
similarly for
$f|_{\pi(\tau(H))}$.  We only need to check that they match up.
For this note that derivatives along $\R^2$ are all zero.
Hence we only need to check one remaining direction.  There we know that
this is along one of the leaves of $H \cup \tau(H)$,
which we know are $C^\infty$ submanifolds.
\end{proof}

To finish the proof of Theorem \ref{flatextthm} we can just apply the
following lemma.  We will use coordinates $(z,w) \in \C^{N-2} \times \C^2$.

\begin{lemma}
Suppose that $U = U_z \times U_w \subset \C^{N-2} \times \C^2$ is a connected
neighbourhood, and 
$\tau(U_w) = U_w$.
$H \subset U$
is a connected Levi-flat $C^\infty$
hypersurface with boundary, with $\partial H \subset \C^{N-2} \times \R^2$.
Then $H \subset \C^{N-2} \times H_w \subset \C^{N-2} \times \C^2$,
where $H_w \subset \C^2$ is a $C^\infty$ Levi-flat hypersurface with boundary
such that $\partial H_w \subset \R^2$.
\end{lemma}

\begin{proof}
We have already seen that the leaves of the foliation induced on $\partial H$
are unions of CR orbits.  Here the CR orbits are just given by
$\{(z,w) \mid w = w^0 \}$ for a fixed $w^0$.
So take one leaf $L$
of the Levi foliation on $H$ extended to the boundary.
It is then easy to see that $L \cap \partial
H$ is equal to (after perhaps extending in the $z$ direction)
to $\C^{N-2} \times A$ for some submanifold
$A \subset \R^2$.

Fix some $p=(z^0,w^0) \in \C^{N-2} \times \R^2$, such that $p \in L$.
It is not hard to see that if we let $L_w :=
\{ w \mid (z^0,w) \in L \}$,
then $L_w \setminus \R^2$ is a codimension 1 complex analytic subvariety
of $V \setminus \R^2$, for some
small neighbourhood $V$ of $w^0$.
Further, $(\C^{N-2} \times L_w)
\cap \partial H = A$, and one component of $L_w$ is
path connected to $A$.  This is because of how $L$ is defined.
If $\tilde{H}$ is any $C^\infty$ submanifold extending $H$ (as noted in the
introduction), then $L$ can be extended to a real $C^\infty$ submanifold
of $\tilde{H}$.  Further, this extension meets
$\C^{N-2} \times \R^2$ transversely in $H$, and all the
derivatives in the $z$ and $\bar{z}$ directions of the defining functions
of $L$ must vanish at $p$,
since $L \cap \partial H = \C^{N-2} \times A$.
Hence, $L_w$ is a submanifold with boundary in some
small neighbourhood of $w^0$.  By dimension, $L$ is then equal to
$\C^{N-2} \times L_w$ in some small neighbourhood
of $p$.  So near some point, $L$ can be defined
by an equation not depending on $z$.  Since $L$ is a connected complex
analytic submanifold, this is true everywhere on $L$.  $H$ is a union
of such $L$ and the lemma follows.
\end{proof}

The uniqueness in Theorem \ref{flatextthm} is obvious in view of the fact
that the extension (near the origin) is given by extension of the leaves of
the Levi foliation and complex submanifolds have unique continuation.

\section{Counterexamples} \label{section:counterex}

In this section we will give examples to show that the assumptions in
Theorems \ref{extthm} and \ref{flatextthm} are indeed optimal.

\begin{example}
It is obvious that Levi-flat hypersurfaces which contain $\R^2 \subset \C^2$
cannot be unique since for example if we have coordinates $(z,w) \in \C^2$,
then both the hypersurfaces $\Im z = 0$ and $\Im w = 0$ contain $\R^2$.
\end{example}

\begin{example} \label{example:smoothH}
We can find a $C^\infty$ Levi-flat hypersurface in $\C^2$ which contains
$\R^2$, but
which is not real analytic (not contained in a real analytic subvariety of
the same dimension).
First let
\begin{equation}
\varphi(x) :=
\begin{cases}
e^{-1/x} & x > 0 ,
\\
0 & x \leq 0 .
\end{cases}
\end{equation}
Then define $H$ by looking at 
\begin{equation} \label{eq:nonrafol}
\rho_t(z,w) := \varphi(t)z^2+t-w .
\end{equation}
On $\R^2$ this defines a $C^\infty$ (but not real analytic)
family of real analytic curves, and it
therefore cannot be induced by a real analytic Levi-flat hypersurface.
We need to show that as
$(z,w)$ range over some neighbourhood of the origin in $\C^2$,
and $t$ ranges over a
small interval, $\rho_t = 0$ defines a Levi-flat
hypersurface.  It suffices to show that it is a submanifold near zero.  It is
automatically Levi-flat since it is given by a 1 parameter family of
complex analytic subvarieties.  First, we check that if $z$, $w$, and $t$ are
kept small, then the complex analytic subvarieties do not intersect for different $t$.
By direct calculation this can be seen to be the case as long as $\abs{z} < 1$.
We look at $\Re \rho_t$ and $\Im \rho_t$,
and notice $\Re \rho_t$ as a function of $(\Re z,\Im z,\Re w,t)$
satisfies the real
analytic implicit
function theorem at 0 and hence we can find a real analytic solution
$t = \alpha(\Re z, \Im z, \Re w)$, then we have a smooth
hypersurface defined by
\begin{equation}
0 = \Im \rho_{\alpha}
= \varphi(\alpha(\Re z, \Im z, \Re w)) \Im (z^2) - \Re w .
\end{equation}
\end{example}

Thus the requirement in Theorem \ref{extthm} that not all
local CR orbits are of codimension 2 in $M$ is necessary.  This is
because the above example extends to $\C^N$ by just letting 
$M = \C^{N-2} \times \R^2$.

\begin{example}
The methods of this paper revolve around extending the Levi foliation
of the hypersurface and thereby extending $H$.  Such methods are 
bound to fail in general when $M$ has a complex tangent and therefore is not
a CR submanifold.  In the following example, we show that even if we can
extend a Levi-flat
hypersurface past a CR singular boundary, the extension need not be unique,
even in the sense of Theorem \ref{flatextthm}.

Let $(z,w) \in \C^2 \times \C$ be our coordinates.
For a fixed $t$, let $H_t$ be a Levi-flat hypersurface defined by
\begin{equation}
\Im w = t\varphi(-\Re w) ,
\end{equation}
where $\varphi$ is as before.  Then define $M$ by
\begin{equation}
\Re w = \abs{z_1}^2 + \abs{z_2}^2 \ \text{ and } \ \Im w = 0 .
\end{equation}
Outside of the origin, $M$ is a CR submanifold, where the codimension in $M$
of the CR orbits must be 1, as $M$ contains no complex analytic
subvarieties.  But then we have a whole family of Levi-flat
hypersurfaces which contain $M$.
\end{example}

\begin{example} \label{example:smooth}
If $M$ would be only $C^\infty$, then no general extension theorem like
Theorem \ref{extthm} nor Theorem \ref{flatextthm} holds.
First, let $\sqrt{\cdot}$ denote the principal branch of the square root,
and note that the function $\xi \mapsto e^{-1/\sqrt{\xi}}$, holomorphic
for $\Re \xi > 0$, can be extended to be $C^\infty$ on $\Re \xi \geq 0$.
Suppose
that in coordinates $(z,w_1,w_2) \in \C^3$ we define a $C^\infty$
Levi-flat hypersurface
with boundary by
\begin{equation}
\Re w_1 \geq \abs{z}^2 \ \text{ and } \ \Re w_2 = \Re e^{-1/\sqrt{w_1}} .
\end{equation}
$M$ is defined similarly by
$\Re w_1 = \abs{z}^2$ and
$\Re w_2 = \Re e^{-1/\sqrt{w_1}}$.
It is easy to check
that $M$ is a generic $C^\infty$ submanifold.  Further,
since $M$ contains no complex analytic subvarieties,
the CR orbits of $M$ can be seen to be of codimension 1 in $M$.
At an interior point, $H$
is given by a vanishing of the real part of a holomorphic function and so $H$
is Levi-flat.

However, $H$ cannot possibly extend across $M$ since that would mean that the
leafs of the
Levi foliation of $H$ would have to extend.  The leaf of $H$ that goes
through the origin is given by
$w_2 = e^{-1/\sqrt{w_1}}$.  Since this subvariety is given as a graph,
if we could possibly extend this complex analytic
subvariety across the origin, we could extend the function
$e^{-1/\sqrt{w_1}}$ across $w_1=0$, and we know this is not possible.
\end{example}

\section{Almost minimal submanifolds} \label{section:almostmin}

We will now prove Theorem \ref{noalmostmin}.  Recall that
a real analytic generic submanifold $M$
is almost minimal at 0 if for every neighbourhood $U$ of 0,
there exists a point $p \in M \cap U$ such that (some representative of)
the local CR orbit at $p$
is not contained in a proper complex analytic subvariety of
$U$.  Let us restate Theorem \ref{noalmostmin} for reader convenience.

\begin{thmnonum} 
Let $M$ be a connected real analytic generic submanifold of
codimension 2 through the origin, which is almost minimal at the origin.
Let $H$ be a connected $C^2$ hypersurface with boundary and $M \subset
\partial H$.  Then $H$ is not Levi-flat.
\end{thmnonum}

Theorem \ref{noalmostmin} is a consequence of the following more general
result.

\begin{lemma} \label{lemma:leavesextend}
Let $M$ be a connected real analytic generic codimension 2 submanifold through
the origin and let $H$ be a connected $C^2$ hypersurface with boundary,
and $M \subset \partial H$.
Suppose that there exists a point on $M$ where the local CR orbits are
of codimension 1 in $M$.  Then there exists some neighbourhood $U$ of the origin
such that the leaves of the Levi foliation of $H^o$ extend to be closed
complex analytic subvarieties of $U$.
\end{lemma}

We will need the following lemma from \cite{Lebl:lfnm}.  Here
$\sX_p$ is the intrinsic complexification of the local CR orbit at $p$,
that is, the smallest germ of a complex analytic submanifold
that contains the local CR orbit at $p$.  
When we say that
$M$ is given in normal coordinates, we mean local holomorphic coordinates
$(z,w) \in \C^{N-2} \times \C^2$, such that $M$ is given near the origin by
\begin{equation}
\begin{split}
& \Im w_1 = \varphi_1(z,\bar{z},\Re w) , \\
& \Im w_2 = \varphi_2(z,\bar{z},\Re w) , \\
\end{split}
\end{equation}
where $\varphi_j (0,\bar{z},s) \equiv \varphi_j (z,0,s) \equiv 0$
for $j = 1,2$.  Thus $M$ is locally a graph over $\C^{N-2} \times \R^2$.

\begin{lemma} \label{Xpnoz}
Given $M \subset U$ in normal coordinates, then there is a small neighbourhood
of the origin $V$ such that for $p \in M \cap V$,
$\sX_p$ contains $\{ (z,w) \in U \mid w = w^0 \}$
as germs at any $(z^0,w^0) \in \sX_p$.
\end{lemma}

We now prove Lemma \ref{lemma:leavesextend} and therefore Theorem
\ref{noalmostmin}.  The method of this proof together with
Theorem \ref{C2extthm} could be used to give a different (but longer)
proof of Theorem \ref{extthm}.

\begin{proof}[Proof of Lemma \ref{lemma:leavesextend}]
We first write $M$ in terms of normal coordinates $(z,w) \in \C^{N-2} \times
\C^2$,
and take $U$
to be the neighbourhood small enough to apply Lemma \ref{Xpnoz}.
(that is $U = V$ in the Lemma).

If the local CR orbits are of codimension 1 in $M$ somewhere on $M$, they are of
codimension 1 outside a proper real analytic subvariety of $M$.  Let $p$
be one of the points where local CR orbits of $M$ are of codimension 1 in $M$.

We note that if $L$ is a leaf of the Levi foliation of $H$ (we extend this
foliation to $M$ as above) such that $p \in L$,
then by Lemma \ref{lemma:uniqueH}, applied in a suitably small
neighbourhood of $p$, we see that as germs $(L,p) \subset \sX_p$. 
Hence we can
extend $L$ to a small neighbourhood of $p$, and it will agree with some
representative of $\sX_p$.
By Lemma \ref{Xpnoz}, we see that near $p$, $L$ is defined by equations
independent of $z$.  Since $L$ is a connected complex submanifold of $U$,
then at each point it is defined by equations independent of $z$.  Hence
there exists a submanifold $\tilde{L}$ of the same dimension,
such that $L \subset \tilde{L}$ and $\tilde{L} = \C^{N-2} \times \tilde{L}_w$
where $\tilde{L}_w$ is a complex hypersurface of $\C^2$.  Now,
if we fix $z = z^0$ and look at $M \cap \{z = z_0\}$,
we see that this is a maximally totally real submanifold of $\C^2$,
and hence locally biholomorphic to $\R^2$.
We can apply the same reasoning as in the proof
of \ref{C2extthm} to apply Schwarz reflection principle
to extend this complex
hypersurface across $\R^2$.
We can therefore assume that $\tilde{L}_w$
is a subvariety of $U \cap \{ z = z_0 \}$ (for a perhaps smaller $U$) and
hence $\tilde{L}$ is a complex analytic subvariety of $U$.
\end{proof}

\section{Almost minimal example} \label{section:ex}

Let $M_\lambda$, $\lambda \in \R$, be the generic, nowhere minimal
submanifold of $\C^3$, with
holomorphic coordinates $(z,w_1,w_2)$ defined by
\begin{equation}
\begin{split}
\bar{w}_1 &= e^{iz\bar{z}} w_1 ,
\\
\bar{w}_2 &= e^{i\lambda z\bar{z}} w_2 .
\end{split}
\end{equation}

When $\lambda$ is irrational, this submanifold is almost minimal at
0, and thus not contained in any Levi-flat real analytic 
subvariety of codimension 1 in $\C^3$,
see \cite{Lebl:lfnm}.  As we will see below, the intrinsic complexification
for a generic point $p = (z^0,w_1^0,w_2^0)$ is given by
\begin{equation}
\{ (z,w) \in \C \times \C^2 \mid 
w_1 = \overline{w_1^0} e^{i(\omega+\theta)},
w_2 = \overline{w_2^0} e^{i\lambda(\omega+\theta)} ) \},
\end{equation}
where $\theta = \arg w_1^0$, and $\omega$ varies over $\C$.  It is not hard
to see that these sets cannot be contained in complex analytic subvarieties
for any neighbourhood at the origin.  To see this note that if we let
$\omega$ vary over $\C$, for any point $p^1 = (z^1,w_1^1,w_2^1)$ in the
set we can (by adding $2\pi$ to $\omega$) get
a dense set of rotations of $w_2^1$ to also be in the set.  This means that
the closure of the set will in general be 5 real dimensional.

When $\lambda = a/b$ is rational, $M_\lambda$
is contained in a Levi-flat subvariety of codimension 1, as
the meromorphic function $w_1^a / w_2^b$ is real valued on
$M_\lambda$.

By Theorem \ref{noalmostmin},
$M_\lambda$ is not a boundary of a $C^2$ Levi-flat hypersurface
for $\lambda$ irrational.
We prove the following theorem to show that it can't be a ``boundary''
of a real analytic Levi-flat subvariety, even if we allow singularities.

\begin{thm} \label{thmex}
Let $\lambda$ be irrational and let $M_\lambda \subset \C^3$ be as above.
Suppose $H$ is a codimension 1 real analytic subvariety
of $D - M_\lambda$, where $D$ is a 
polydisc in $\C^3$ centered at the origin.
Suppose that there exists a point $p \in M_\lambda$,
and a connected $C^2$ hypersurface $N$ with boundary,
such that $p \in \partial N \subset M_\lambda$, and $N^o \subset H$.
Then $H$ is not Levi-flat.

In fact, if $H$ is irreducible, then $H$ is not Levi-flat
at any nonsingular point of top dimension.
\end{thm}

\begin{proof}
Assume for contradiction that $H$ is Levi-flat.  In particular this means
that $N^o \cap H^*$ is Levi-flat, where $H^*$ are the nonsingular points of
hypersurface dimension.  Thus $N^o$ is Levi-flat on an open dense set.
Since being Levi-flat means a certain $C^1$ 1-form is integrable, then 
it is integrable on all of $N^o$ by continuity.

Pick a point $p = (z, w_1 , w_2)$ on $M =
M_\lambda$, such that $p \in N$, and the
local CR orbits of $M$ are of codimension 1 in $M$
in a neighbourhood $U \subset M$
of $p$.

If we take $q \in U$, and $\sX_q$ is (some representative of the germ of)
intrinsic complexification of the
local CR orbit, then $M \cap \sX_q$
is a hypersurface in $\sX_q$ and hence divides $\sX_q$ into two connected
sets (we can pick a representative of $\sX_q$ small enough).
Hence we can write $\sX_q$ as a disjoint union of three connected sets as
follows:
\begin{equation}
\sX_q = \sX_q^+ \cup (M \cap \sX_q) \cup \sX_q^- .
\end{equation}
By Lemma \ref{lemma:uniqueH}, we see that either
$\sX_q^+ \subset N^o$ or $\sX_q^- \subset N^o$.  So suppose
$\sX_q^+ \subset N^o \subset H$.

Now we will find a parametrization of $\sX_q$ and hence of $\sX_q^+$.
We will construct this parametrization of $\sX_q$ by the use of Segre sets.
We can compute the third Segre set at $q = (z^0,w_1^0,w_2^0)$,
where $w_1^0 \not= 0$ and $w_2^0 \not= 0$,
by the following mapping (see \cite{BER:book})
\begin{equation}
\varphi(t_1,t_2,t_3) :=
(t_3,
\overline{w_1^0} e^{i(t_3 t_2 - t_2 t_1 + t_1 \overline{z^0})},
\overline{w_2^0} e^{i\lambda(t_3 t_2 - t_2 t_1 + t_1 \overline{z^0})} ) .
\end{equation}
That is, the image of this mapping agrees with $\sX_q$
as germs at $p$.
We must be careful to stay
within the polydisc $D \subset \C^3$ in the image.
So let us suppose that
$M$ is only defined in $D$.

Let $\theta := \arg w_1^0$.  On $M$, $\arg w_1^0 = \frac{1}{\lambda} \arg w_2^0 = \theta$.
Changing variables by precomposing with $(\xi,\omega) \mapsto
(0,\frac{\omega+\theta}{\xi},\xi)$ we get the map:
\begin{equation}
\tilde{\varphi}(\xi,\omega) :=
(\xi,
\overline{w_1^0} e^{i(\omega+\theta)},
\overline{w_2^0} e^{i\lambda(\omega+\theta)} ) .
\end{equation}
The image of this map is on $M$ when
\begin{equation}
\begin{split}
\overline{w_1^0} e^{i(\omega+\theta)}
& =
e^{i\abs{\xi}^2}
w_1^0 e^{-i(\overline{\omega}+\theta)},
\\
\overline{w_2^0} e^{i\lambda(\omega+\theta)}
& =
e^{i\lambda\abs{\xi}^2}
w_2^0 e^{-i\lambda(\overline{\omega}+\theta)} .
\end{split}
\end{equation}
That is, the pullback of the CR orbit at $q$ by $\tilde{\varphi}$ is
\begin{equation} \label{eq:Sdef}
\Re\omega =
\frac{1}{2}\abs{\xi}^2 .
\end{equation}
Let $S$ be hypersurface
in the parameter space $(\xi,\omega)$ defined by \eqref{eq:Sdef}.
So if, in the parameter space $(\xi,\omega)$, we stay on one or the other
side of $S$,
we are parametrizing either $\sX_q^+$ or $\sX_q^-$.

Let us also vary $w_1^0$ and $w_2^0$, while keeping $q=(z^0,w_1^0,w_2^0)$
within $M$.
That is let $w_1^0 = re^{i\theta}$ and
$w_2^0 = se^{i\lambda\theta}$, and now let $r$ and $s$ vary.
We define a map $\psi$
by adding the parameters $r$ and $s$ to $\tilde{\varphi}$
\begin{equation} \label{eq:psidef}
\psi(\xi,\omega,r,s) :=
(\xi,
r e^{-i\theta}
e^{i(\omega+\theta)},
s e^{-i\lambda\theta}
e^{i\lambda(\omega+\theta)} ) =
(\xi,
r
e^{i\omega},
s
e^{i\lambda\omega} ) .
\end{equation}

As $r$ and $s$ vary over a small interval
and $\xi$ and $\omega$ vary over some small connected
open set, such that the image of $\psi$ never leaves $D$,
and further, such that $\xi$ and $\omega$ stay on one
side of $S$, we get a parametrization of an open part of $H$.
This is because as we vary $r$ and $s$, we vary $q$, and then as
we vary $\xi$ and $\omega$, we parametrize $\sX_q^+$ (as long as
$\xi$ and $\omega$ stay on one side of $S$).

We will make the parametrization an immersion
by restricting
$\omega$ to be real.  Then for a small open set
$V \subset \C \times \R^3$, $\psi|_V$ is an immersion.
We pick this $V$ such that $\psi(V) \subset H$.
Now pick any connected open $V' \subset \C \times \R^3$,
such that $V \subset V'$ and for all
$(\xi,\omega,r,s) \in V'$ we have
$\omega > \frac{1}{2}\abs{\xi}^2$ (or $\omega <
\frac{1}{2}\abs{\xi}^2$) and $r,s \in (0,\epsilon)$ (where $\epsilon$ 
is the radius of $D$).  It is clear that $\psi(V') \subset D \setminus M$.
Further, $\psi(V') \subset H$, since $H$
is a subvariety of $D \setminus M$ and $V'$ is connected and if we pull back
$H$ by $\psi|_{V'}$ we must get a subvariety of $V'$ which contains $V$.


Note that we can pick $V'$ such that it contains all
$\omega \in (0,\infty)$ (or in $(-\infty,0)$).  Without loss of generality
suppose we can let $\omega$ go to plus infinity and still stay within $H$

Now we will show that
$H$ must be dense (in $\C^3$) near some point
not on $M$, but arbitrarily close to 0.
Let $\xi$ vary in some small open set and let $r$ and $s$ vary in some small
open interval.  For a bounded interval of $\omega$ we will parametrize a 5
dimensional set.  Now we can start adding $2\pi$ to $\omega$ and we add a
dense set of rotations to the third component in \eqref{eq:psidef},
without changing the first two.  Thus the image of $\psi$ must be dense near
some point and this contradicts $H$ being a subvariety of codimension 1.


Now suppose that $H$ is irreducible.
Let $H^*$ be the nonsingular points of top dimension.

In \cite{burnsgong:flat} (Lemmas 2.1 and 2.2),
Burns and Gong prove the following.
Let $K$ be a subvariety of codimension 1 ($0 \in K$)
defined by $r(z,\bar{z}) = 0$, for $r$
an irreducible
real analytic real valued function.  Then for some small neighbourhood
$U$ of 0, $r$ complexifies (the Taylor series $r(z,w)$ converges for
$z \in U$, $\bar{w} \in U$) and is irreducible as a holomorphic function.
Further, if $K^* \cap U$ is Levi-flat at a single point, then $K^* \cap U$
is Levi-flat at all points.

We can use this to show that if $H^*$ is Levi-flat at one point
and $H$ is irreducible in $D \setminus M$, then $H^*$ is Levi-flat
at all points and hence $H$ is Levi-flat by our definition.  By the above
result we can find a collection of open neighbourhoods
$U_j$ and for each $U_j$ we find irreducible branches $A_{j1},\ldots,A_{jn}$
of $H$ in $U_j$, and assume that each $A_{jk} \subset U_j$ satisfies the
above property.  Now take $H'$ be a union of those $A_{jk}$ such that
$A_{11} \subset H'$ and if $A_{jk} \subset H'$ and
$A_{\ell m} \cap A_{jk}$ is of codimension 1, then $A_{\ell m} \subset H'$.
It is clear that $H'$ is a subvariety of $D \setminus M$ and since $H$
is irreducible then $H' = H$.  It is clear that all the $A_{jk}$ are
Levi-flat if and only if $A_{11}$ is Levi-flat, and we are done.
\end{proof}


\section{Subanalytic hypersurfaces} \label{section:subanal}

If we allow subanalytic hypersurfaces (see
\cite{BM:semisub}), then we have the following result.

\begin{thm} \label{thmsub}
Let $M$ be a real analytic, codimension 2, generic submanifold that
is nowhere minimal.  Then there exists a subanalytic hypersurface $H$,
which is Levi-flat at nonsingular points, such that $M \subset H$.
Further, if $H^*$ are the nonsingular points of top dimension of $H$, then $M \cap H^*$ is dense
in $M$.
\end{thm}

\begin{proof}
If all CR orbits of $M$ are of codimension 2 in $M$, this is trivial.
Otherwise, intersect with a small ball
around any point in which normal coordinates $(z,w)$ are defined.  Then
take the projection $\pi_w$ onto the $w$ factor.
$\pi_w(M)$ is a subanalytic
hypersurface in general.
Apply Lemma \ref{Xpnoz} to see that all the $\sX_p$ are product sets,
and $\pi_w(\sX_p)$ is contained in $\pi_w(M)$.
If $\pi_w(M)$ is of codimension 2, $M$ had CR orbits of only codimension 2 in
$M$.
If $\pi_w(M)$ is of codimension 0, then $M$ must have been minimal.
Hence $\pi_w(M)$ must have been a subanalytic hypersurface foliated
by complex analytic subvarieties (the projections of the CR orbits),
since $M$ is nowhere minimal.  Thus
$\pi_w(M)$ is the subanalytic hypersurface we are looking for.
See \cite{Lebl:lfnm} for more details of this method.
\end{proof}

Note that we must intersect with a small ball first, else the image of
the projection need not be subanalytic.  The submanifold $M_\lambda$
for $\lambda$ irrational from
\S \ref{section:ex}, when projected onto
the $w$ factor without restricting the $z$ to be bounded, will be a
dense set in $\C^2$ which is not subanalytic.  On the other hand,
if we intersect $M_\lambda$ with the set $\abs{z} \leq 1$,
and we look at $\pi_w(M_\lambda)$,
we get the following
subanalytic hypersurface:
\begin{equation}
\left\{ w \in \C^2 ~\Big|~
\arg w_1 = t, ~ \arg w_2 = \lambda t, ~ -\frac{1}{2} \leq t \leq 0 \right\} .
\end{equation}

Note that $H^*$ (the nonsingular points of top dimension of $H$)
is again a subanalytic set and hence a locally finite union of
real analytic submanifolds.
Thus we have that a nowhere minimal
$M$ is contained in the closure of 
a locally finite union of real analytic Levi-flat hypersurfaces.
If the points where
the CR foliation of $M$ is of codimension 1 are connected, then we need only
take one hypersurface.  However, $H$ need not have smooth boundary
nor does the boundary need to be equal to $M$ if it does.
Thus we cannot apply Theorem \ref{extthm}.

If we allow hypersurfaces with singularities all the way up to $M$
in the sense of \cite{DTZ:boundfull}, then the above result suggests
that, at least locally, any nowhere minimal submanifold could conceivably
bound such a singular hypersurface.

\renewcommand\MR[1]{\relax\ifhmode\unskip\spacefactor3000 \space\fi
  \def\@tempa##1:##2:##3\@nil{%
    \ifx @##2\@empty##1\else\textbf{##1:}##2\fi}%
  \href{http://www.ams.org/mathscinet-getitem?mr=#1}{MR \@tempa#1:@:\@nil}}


\begin{bibdiv}
\begin{biblist}

\bib{BER:book}{book}{
      author={Baouendi, M.~Salah},
      author={Ebenfelt, Peter},
      author={Rothschild, Linda~Preiss},
       title={Real submanifolds in complex space and their mappings},
      series={Princeton Mathematical Series},
   publisher={Princeton University Press},
     address={Princeton, NJ},
        date={1999},
      volume={47},
        ISBN={0-691-00498-6},
      review={\MR{2000b:32066}},
}

\bib{BG:envhol}{article}{
      author={Bedford, Eric},
      author={Gaveau, Bernard},
       title={Envelopes of holomorphy of certain {$2$}-spheres in {${\bf
  C}\sp{2}$}},
        date={1983},
        ISSN={0002-9327},
     journal={Amer. J. Math.},
      volume={105},
      number={4},
       pages={975\ndash 1009},
      review={\MR{84k:32016}},
}

\bib{BM:semisub}{article}{
      author={Bierstone, Edward},
      author={Milman, Pierre~D.},
       title={Semianalytic and subanalytic sets},
        date={1988},
        ISSN={0073-8301},
     journal={Inst. Hautes \'Etudes Sci. Publ. Math.},
      number={67},
       pages={5\ndash 42},
      review={\MR{89k:32011}},
}

\bib{Bishop:diffman}{article}{
      author={Bishop, Errett},
       title={Differentiable manifolds in complex {E}uclidean space},
        date={1965},
        ISSN={0012-7094},
     journal={Duke Math. J.},
      volume={32},
       pages={1\ndash 21},
      review={\MR{34:369}},
}

\bib{Boggess:CR}{book}{
      author={Boggess, Albert},
       title={C{R} manifolds and the tangential {C}auchy-{R}iemann complex},
      series={Studies in Advanced Mathematics},
   publisher={CRC Press},
     address={Boca Raton, FL},
        date={1991},
        ISBN={0-8493-7152-X},
      review={\MR{94e:32035}},
}

\bib{burnsgong:flat}{article}{
      author={Burns, Daniel},
      author={Gong, Xianghong},
       title={Singular {L}evi-flat real analytic hypersurfaces},
        date={1999},
        ISSN={0002-9327},
     journal={Amer. J. Math.},
      volume={121},
      number={1},
       pages={23\ndash 53},
      review={\MR{2000j:32062}},
}

\bib{DAngelo:CR}{book}{
      author={D'Angelo, John~P.},
       title={Several complex variables and the geometry of real
  hypersurfaces},
      series={Studies in Advanced Mathematics},
   publisher={CRC Press},
     address={Boca Raton, FL},
        date={1993},
        ISBN={0-8493-8272-6},
      review={\MR{94i:32022}},
}

\bib{DTZ:boundfull}{unpublished}{
      author={Dolbeault, Pierre},
      author={Tomassini, Giuseppe},
      author={Zaitsev, Dmitri},
       title={On boundaries of {L}evi-flat hypersurfaces in
  {${\mathbb{C}}^n$}},
        note={preprint},
}

\bib{Flanders:dforms}{book}{
      author={Flanders, Harley},
       title={Differential forms with applications to the physical sciences},
     edition={Second},
      series={Dover Books on Advanced Mathematics},
   publisher={Dover Publications Inc.},
     address={New York},
        date={1989},
        ISBN={0-486-66169-5},
      review={\MR{90k:53001}},
}

\bib{Lebl:lfnm}{article}{
      author={Lebl, Ji{\v r}\'i},
       title={Nowhere minimal {C}{R} submanifolds and {L}evi-flat
  hypersurfaces},
     journal={J. Geom. Anal.},
      pages = {321--342},
     volume = {17},
       year = {2007},
     number = {2},
        note={arXiv:
  \href{http://www.arxiv.org/abs/math.CV/0606141}{math.CV/0606141}},
}

\bib{Malgrange}{book}{
      author={Malgrange, B.},
       title={Ideals of differentiable functions},
      series={Tata Institute of Fundamental Research Studies in Mathematics,
  No. 3},
   publisher={Tata Institute of Fundamental Research},
     address={Bombay},
        date={1967},
      review={\MR{35:3446}},
}

\bib{MW:normal}{article}{
      author={Moser, J{\"u}rgen~K.},
      author={Webster, Sidney~M.},
       title={Normal forms for real surfaces in {${\bf C}\sp{2}$} near complex
  tangents and hyperbolic surface transformations},
        date={1983},
        ISSN={0001-5962},
     journal={Acta Math.},
      volume={150},
      number={3-4},
       pages={255\ndash 296},
      review={\MR{85c:32034}},
}

\bib{Nagano}{article}{
      author={Nagano, Tadashi},
       title={Linear differential systems with singularities and an application
  to transitive {L}ie algebras},
        date={1966},
     journal={J. Math. Soc. Japan},
      volume={18},
       pages={398\ndash 404},
      review={\MR{33:8005}},
}

\bib{StraubeSucheston:fol}{article}{
      author={Straube, Emil~J.},
      author={Sucheston, Marcel~K.},
       title={Levi foliations in pseudoconvex boundaries and vector fields that
  commute approximately with {$\overline\partial$}},
        date={2003},
        ISSN={0002-9947},
     journal={Trans. Amer. Math. Soc.},
      volume={355},
      number={1},
       pages={143\ndash 154 (electronic)},
      review={\MR{2003h:32058}},
}

\bib{Sussmann}{article}{
      author={Sussmann, H{\'e}ctor~J.},
       title={Orbits of families of vector fields and integrability of
  distributions},
        date={1973},
        ISSN={0002-9947},
     journal={Trans. Amer. Math. Soc.},
      volume={180},
       pages={171\ndash 188},
      review={\MR{47:9666}},
}

\bib{Tumanov}{article}{
      author={Tumanov, A.~E.},
       title={Extension of {CR}-functions into a wedge from a manifold of
  finite type},
        date={1988},
        ISSN={0368-8666},
     journal={Mat. Sb. (N.S.)},
      volume={136(178)},
      number={1},
       pages={128\ndash 139},
      review={\MR{89m:32027}},
}

\end{biblist}
\end{bibdiv}
\end{document}